\documentclass{amsart}
\usepackage{amsmath,amsthm,amssymb}
\usepackage{graphicx}
\newtheorem{theorem}{Theorem}[section]
\newtheorem{proposition}[theorem]{Proposition}
\newtheorem{lemma}[theorem]{Lemma}
\newtheorem{corollary}[theorem]{Corollary}
\theoremstyle{definition}
\newtheorem{definition}[theorem]{Definition}
\newtheorem{remark}[theorem]{Remark}

\newtheorem{ack}{Acknowledgement}

\title{Small exotic Stein manifolds}
\author{Selman Akbulut and Kouichi Yasui}
\thanks{The first author is partially supported by NSF, and the second author is partially supported by JSPS Research Fellowships for Young Scientists.}
\date{February 23, 2009}
\address{Department~of~Mathematics, Michigan State University, E.Lansing, MI, 48824, USA}
\email{akbulut@math.msu.edu}

\address{Department~of~Mathematics, Graduate~School~of~Science, Osaka~University, Toyonaka, Osaka 560-0043, Japan}
\email{kyasui@cr.math.sci.osaka-u.ac.jp}

\begin{document}

\begin{abstract}
It is known that the only Stein filling of the standard contact structure on $S^3$ is $B^4$. In this paper, we construct pairs of homeomorphic but not diffeomorphic simply connected compact Stein $4$-manifolds, for any Betti number $b_{2}\geq 1$; we do this by enlarging corks and plugs. 
\end{abstract}

\maketitle

\section{Introduction} A properly imbedded complex submanifold of an affine space $X\subset {\bf C}^N$ is called  a {\it Stein manifold}. A compact smooth submanifold $M\subset X$ is called a {\it Compact Stein manifold} if it is cut out from $X$ by $f\leq c$, where $f:X\to {\bf R}$ is a strictly pluri-subharmonic (proper) Morse function, and $c$ is a regular value.  In particular $M$ is a symplectic manifold with convex boundary and the symplectic form $\omega =\frac{1}{2}\partial\bar{\partial}f$. The form $\omega$ induces a contact structure $\xi$ on the boundary $\partial M$. We call $(M,\omega)$ a {\it Stein filling} of the boundary contact manifold $(\partial M, \xi)$. Stein manifolds have been a useful tools for studying exotic smooth structures on $4$-manifolds, since  smooth $4$-manifolds can be decomposed into codimension zero Stein pieces (e.g \cite{AM1}, ~\cite{AY}).

\vspace{.05in}

 In \cite{E1},~\cite{E2} Eliashberg characterized the topology of Stein manifolds and proved that $S^3$ with the standard contact structure has a unique Stein filling $B^4$ (for more uniqueness results see \cite{OS1}). However, Ozbagci-Stipsicz~\cite{OS2} and Smith~\cite{S} constructed infinitely many different Stein fillings of some contact $3$-manifolds up to diffeomorphism. Their Stein fillings are distinguished by their first homology groups, so in particular they are not homeomorphic to each other. It is thus interesting to find pairs of homeomorphic but not diffeomorphic compact Stein $4$-manifolds. Recently,  by using Fintushel-Stern's knot surgery, Akhmedov-Etnyre-Mark-Smith~\cite{AEMS} constructed infinitely many simply connected compact Stein $4$-manifolds that are mutually homeomorphic but not diffeomorphic, and the induced contact structures on their boundaries are isomorphic. Here we give various  examples of  small simply connected compact Stein manifold pairs, that are homeomorphic but not diffeomorphic to each other, and describe their handlebodies concretely. In the light of the Eliashberg's theorem it is interesting to seek such exotic pairs with the smallest second Betti number $b_{2}$. Our examples provide the case $b_2$=1. Note that the previously constructed exotic contractible manifold of \cite{A1} is only exotic rel boundary, and the exotic manifold pairs of \cite{A2} have $b_{2}=1$ but are not  both Stein.
 
\vspace{.1in}

Corks and plugs are fundamental objects detecting exotic smooth structures on $4$-manifolds \cite{AY}. In~\cite{A1} the first example of a cork was found. In~\cite{A2} an exotic pair of simply connected $4$-manifolds with boundaries was constructed, by enlarging a cork. In~\cite{AY} examples of infinitely many corks and plugs were found. In this paper, we construct examples of exotic pairs of compact Stein $4$-manifolds, by enlarging corks and plugs. We also give handlebody pictures of these Stein manifolds. 
\begin{theorem}\label{th:1.1}
For every second Betti number $b_{2}\geq 1$, 
there exist pairs of simply connected compact Stein $4$-manifolds which are homeomorphic but not diffeomorphic to each other. 
\end{theorem}
\begin{theorem}\label{th:1.2}
For every $n\geq 0$, there exist pairs of simply connected compact Stein $4$-manifolds with the following properties:\\
$(\textnormal{i})$ They are homeomorphic but not diffeomorphic to each other;\\
$(\textnormal{ii})$ The first homology groups of their boundaries are  $\mathbf{Z}/n\mathbf{Z}$. 
\end{theorem}


\begin{theorem}\label{th:1.3}
There exist pairs of simply connected compact Stein $4$-manifolds with the following properties:\\
$(\textnormal{i})$ Their boundaries are diffeomorphic;\\
$(\textnormal{ii})$ Their integer coefficient homology groups are isomorphic;\\
$(\textnormal{iii})$ Their intersection forms are not isomorphic. In particular, they are not homeomorphic. 
\end{theorem}

Let $X$ be a simply connected compact smooth $4$-manifold with $H_2(X;\mathbf{Z})\cong \mathbf{Z}$. Define the genus $G(X)$ of $X$ as the minimal number of the genera of surfaces which represent a generator of $H_2(X;\mathbf{Z})$. In this paper, by enlarging corks we also construct the examples of exotic pairs below. 
These examples show that the distance of genera of two homeomorphic (Stein) $4$-manifolds can be arbitraly large. 
\begin{theorem}\label{th:1.4}
For every $n\geq 0$, there exist simply connected compact Stein $4$-manifolds $X_n$ and $Y_n$ with the following properties:\\
$(\textnormal{i})$ They are homeomorphic but not diffeomorphic to each other;\\
$(\textnormal{ii})$ $H_2(X_n;\mathbf{Z})\cong H_2(Y_n;\mathbf{Z})\cong \mathbf{Z}$;\\
$(\textnormal{iii})$ $G(X_n)-G(Y_n)\geq n$. 
\end{theorem}
 
\begin{ack}The second author would like to thank his adviser Hisaaki Endo for his heartfelt encouragement. This work was partially done during the second author's stay at Michigan State University. The second author is greatful for their hospitality. Finally, the authors would like to thank the referee for his/her useful comments. 
\end{ack}
\section{Corks and Plugs}\label{sec:cork}
In this section, we briefly recall corks and plugs. For details, see \cite{AY}. For basics of handlebody pictures, see \cite{GS}. 
\begin{definition}
Let $C$ be a compact Stein $4$-manifold with boundary, and  let $\tau: \partial C\to \partial C$ be an involution on the boundary. 
We call $(C, \tau)$ a \textit{Cork} if $\tau$ extends to a self-homeomorphism of $C$, but cannot extend to any self-diffeomorphism of $C$. 
A cork $(C, \tau)$ is called a cork of a smooth $4$-manifold $X$, if  $C\subset X$ and $X$ changes its diffeomorphism type when removing $C$ and regluing it via $\tau$. Note that this operation does not change the homeomorphism type of $X$.
\end{definition}
\begin{definition}
Let $P$ be a compact Stein $4$-manifold with boundary, and let $\tau: \partial P\to \partial P$ be an involution on the boundary, which cannot extend to any self-homeomorphism of $P$. We call $(P, \tau)$ a \textit{Plug} of $X$,  if  $P\subset X$ and $X$ keeps its homeomorphism type and changes its diffeomorphism type when removing $P$ and gluing it via $\tau$. 
We call $(P, \tau)$ a \textit{Plug} if there exists a smooth $4$-manifold $X$ such that $(P, \tau)$ is a plug of $X$. 
\end{definition}

\begin{definition} Let $W_n$ and $W_{m,n}$ be the smooth 4-manifolds in Figure~$\ref{fig1}$. 
Let $f_n:\partial W_n\to \partial W_n$ and $f_{m,n}:\partial W_{m,n}\to \partial W_{m,n}$ be the obvious involutions obtained by first surgering $S^1\times B^3$ to $B^2\times S^2$ in the interiors of $W_n$ and $W_{m,n}$, then surgering the other imbedded $B^2\times S^2$ back to $S^1\times B^3$ (i.e. replacing the dots in Figure ~$\ref{fig1}$), respectively. Note that the diagrams of $W_n$ and $W_{m,n}$ are symmetric links.
\begin{figure}[ht!]
\begin{center}
\includegraphics[width=3.5in]{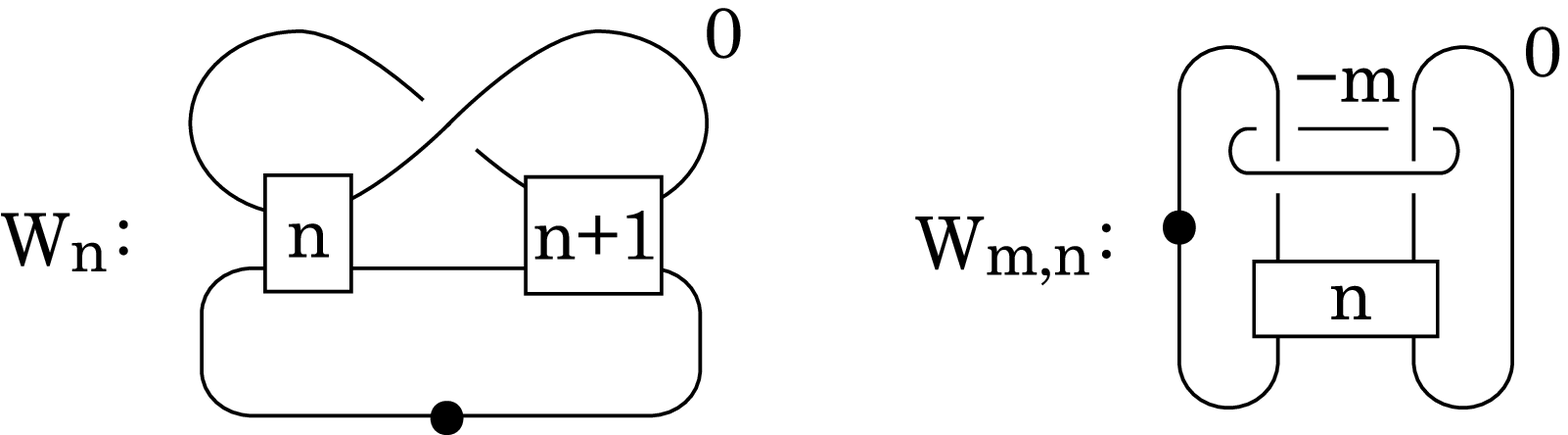}
\caption{}
\label{fig1}
\end{center}
\end{figure}
\end{definition}
\begin{theorem}[{\cite{A1}, \cite{AY}}]\label{th:cork}
$(1)$ For $n\geq 1$, the pair $(W_n, f_n)$ is a cork. \\
$(2)$ For $m\geq 1$ and $n\geq 2$, the pair $(W_{m,n}, f_{m,n})$ is a plug.
\end{theorem}
\section{Proofs}
\subsection{Enlarging corks}In this subsection, we construct exotic smooth pairs of compact Stein $4$-manifolds by enlarging corks $W_n$.

\begin{definition}
Let $C_1(m,n,p,q)$ be the simply connected compact $4$-manifold in Figure~$\ref{fig2}$. Let $C_2(m,n,p,q)$ be the simply connected compact $4$-manifold in Figure~$\ref{fig3}$. Here each $m$-framed knot of Figure~$\ref{fig2}$ and $\ref{fig3}$ is the $(p,p-1)$-torus knot. 
\begin{figure}[ht!]
\begin{center}
\includegraphics[width=4.3in]{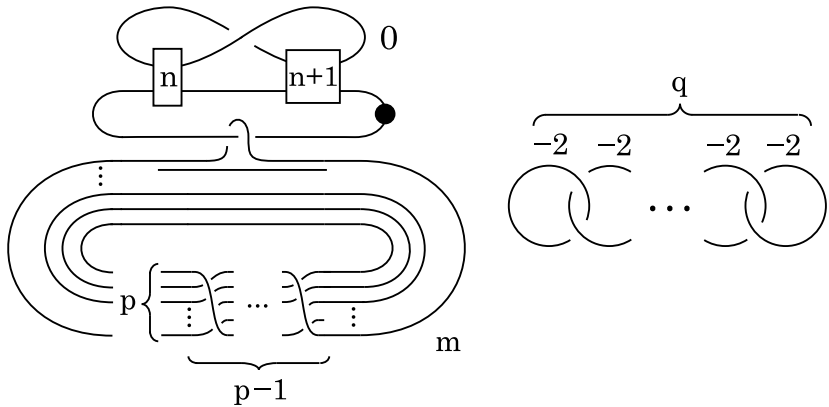}
\caption{$C_1(m,n,p,q)$}
\label{fig2}
\end{center}
\end{figure}
\begin{figure}[ht!]
\begin{center}
\includegraphics[width=4.3in]{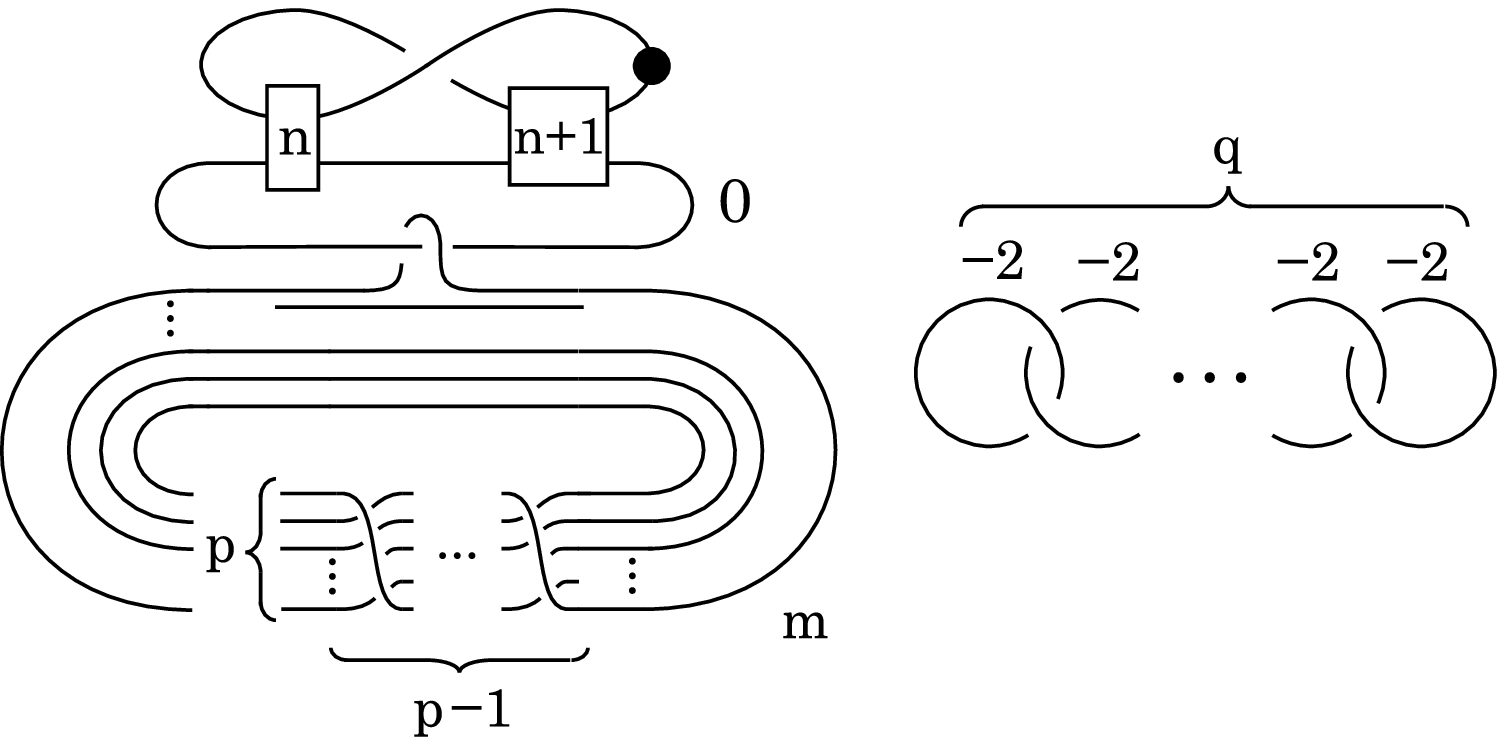}
\caption{$C_2(m,n,p,q)$}
\label{fig3}
\end{center}
\end{figure}
\end{definition}
\begin{lemma}\label{lem:basic}
$(1)$ $C_1(m,n,p,q)$ and $C_2(m,n,p,q)$ are homeomorphic. \\
$(2)$ $C_1(m,n,p,q)$ and $C_2(m,n,p,q)$ are compact Stein $4$-manifolds for $n\geq 1$, $q\geq 0$ and $m\leq p^2-3p$. \\
$(3)$ $H_1(\partial C_1(m,n,p,0);\mathbf{Z})\cong H_1(\partial C_2(m,n,p,0);\mathbf{Z})\cong \mathbf{Z}/m\mathbf{Z}$. \\
$(4)$ $H_2(C_1(m,n,p,q);\mathbf{Z})\cong H_2(C_2(m,n,p,q);\mathbf{Z})\cong \oplus_{q+1} \mathbf{Z}$. 
\end{lemma}
\begin{proof}
(1) $C_2(m,n,p,q)$ is obviously obtained from $C_1(m,n,p,q)$ by removing the contractible $4$-manifold $W_n$ and regluing it via the involution $f_n$. The claim now follows from Freedman's theorem. 

(2) By isotoping the diagrams of $C_1(m,n,p,q)$ and $C_2(m,n,p,q)$ into Legendrian diagrams (and changing the $1$-handle notation) as in Figure~\ref{fig4} and~\ref{fig5}, respectively, we can easily see that these manifolds are Stein, by checking the Eliashberg's criterium: The framing of each 2-handle is less than the Thurston-Bennequin number.

(3) The claim follows from the fact that the intersection forms of $C_1(m,n,p,0)$ and $C_2(m,n,p,0)$ are isomorphic to $\langle m \rangle$. 

(4) This clearly follows from the definition. 
\end{proof}
\begin{figure}[ht!]
\begin{center}
\includegraphics[width=4.8in]{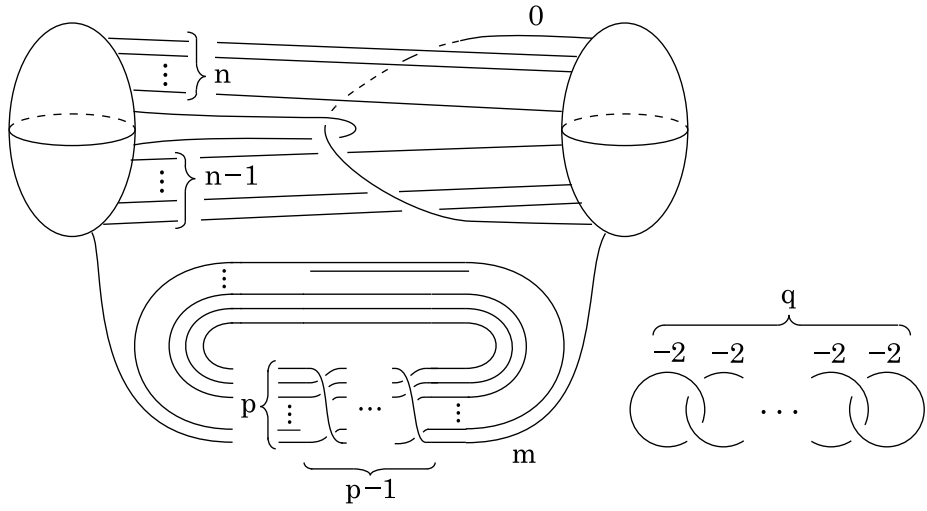}
\caption{$C_1(m,n,p,q)$}
\label{fig4}
\end{center}
\end{figure}
\begin{figure}[ht!]
\begin{center}
\includegraphics[width=4.6in]{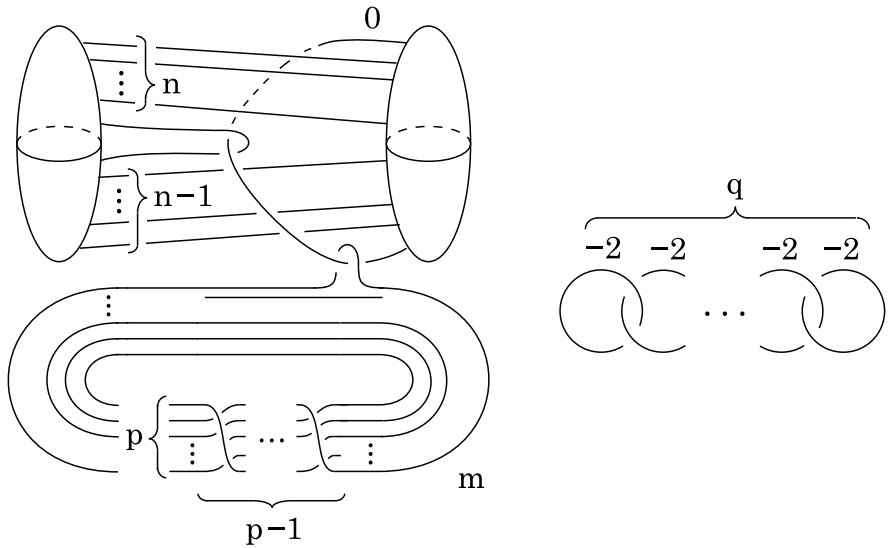}
\caption{$C_2(m,n,p,q)$}
\label{fig5}
\end{center}
\end{figure}
 We first prove the theorem below. This theorem together with Lemma~\ref{lem:basic} clearly gives Theorem~\ref{th:1.1} and \ref{th:1.2}. 
\begin{theorem}\label{th:exotic by cork}
$(1)$ $C_1(m,n,p,0)$ and $C_2(m,n,p,0)$ are homeomorphic but not diffeomorphic to each other, for $n\geq 4$, $p\geq 1$ and $m\leq p^2-3p+1$. \medskip \\
$(2)$ $C_1(m,n,p,0)$ and $C_2(m,n,p,0)$ are homeomorphic but not diffeomorphic to each other, for $1\leq n\leq 3$, $p\geq 3$ and $m\leq p^2-3p+1$. \medskip \\
$(3)$ $C_1(m,n,p,q)$ and $C_2(m,n,p,q)$ are homeomorphic but not diffeomorphic to each other, for $n\geq 1$, $p,q\geq 1$ and $0\leq m\leq p^2-3p+1$. 
\end{theorem}
Let $E(n)$ be the relatively minimal simply connected elliptic surface with Euler characteristic $12n$. We begin with the theorem below. 
\begin{theorem}[Gompf-Stipsicz~\cite{GS}]\label{th:GS}
For $n\geq 1$, the elliptic surface $E(n)$ has the handle decomposition in Figure~$\ref{fig6}$. The obvious cusp neighborhood $($i.e. the dotted circle, $-1$-framed meridian of the dotted circle, and the left most $0$-framed unknot$)$ is isotopic to the regular neighborhood of a cusp fiber of $E(n)$. 
\begin{figure}[ht!]
\begin{center}
\includegraphics[width=3.25in]{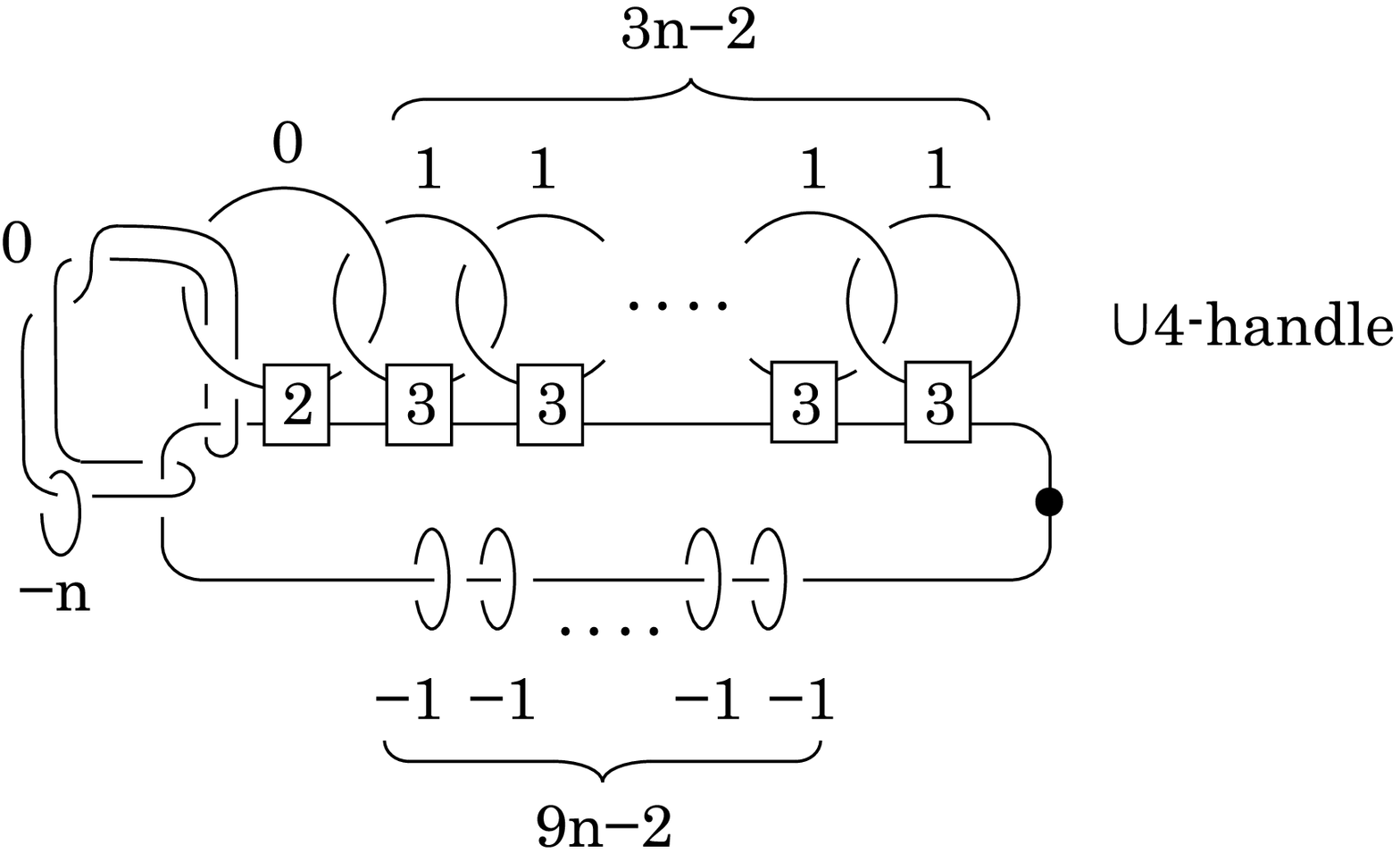}
\caption{$E(n)$}
\label{fig6}
\end{center}
\end{figure}
\end{theorem}
\begin{corollary}\label{cor:GS}
For $n\geq 1$, the elliptic surface $E(n)$ has a handle decomposition as in Figure~$\ref{fig7}$. The obvious cusp neighborhood is isotopic to the regular neighborhood of a cusp fiber of $E(n)$. 
\begin{figure}[ht!]
\vspace{1.0\baselineskip}
\begin{center}
\includegraphics[width=3.55in]{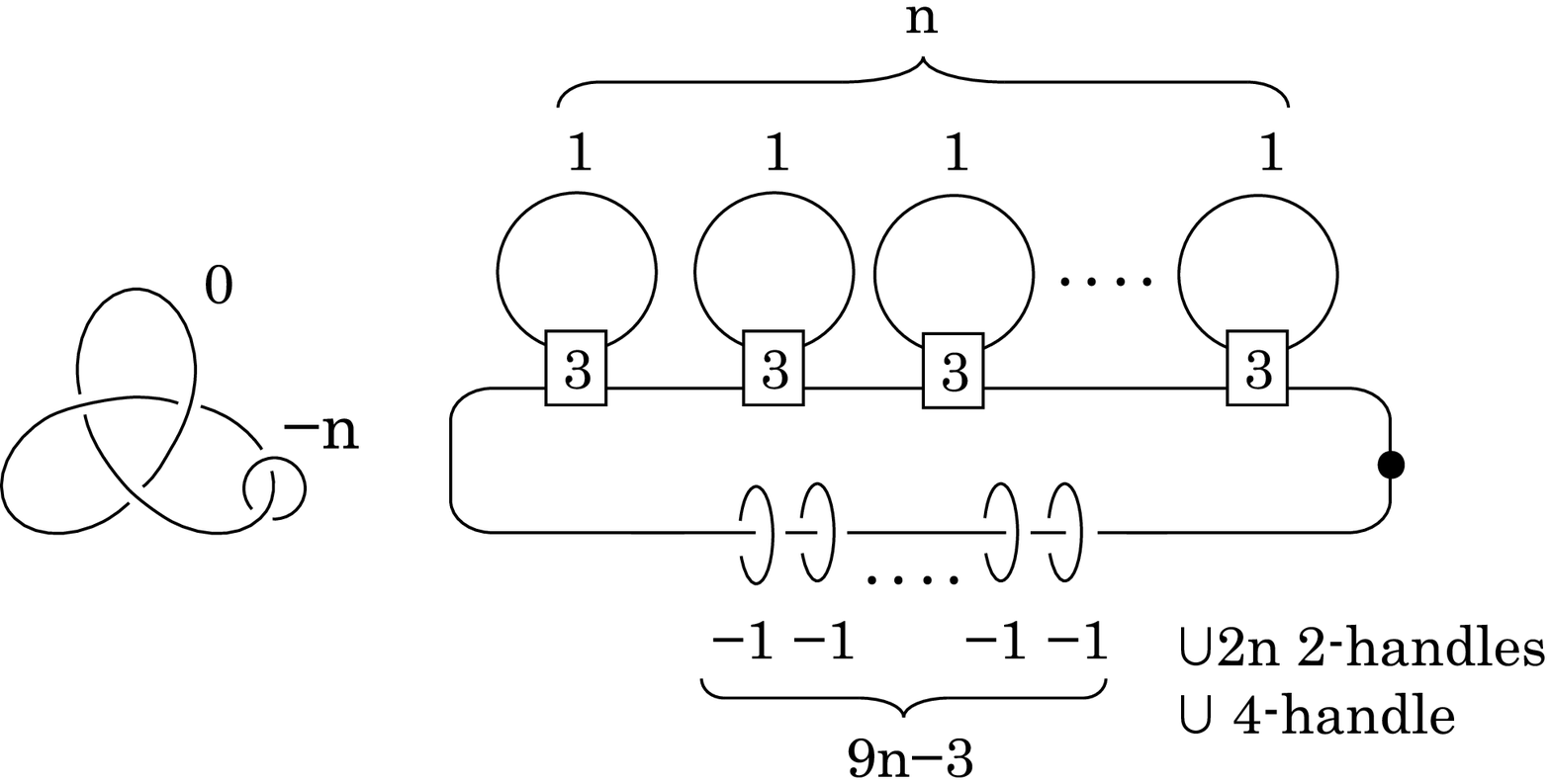}
\caption{$E(n)$}
\label{fig7}
\end{center}
\end{figure}
\end{corollary}
\begin{proof}
In Figure~$\ref{fig6}$, pull off the leftmost 0-framed unknot from the dotted circle by sliding over $-1$-framed knot. 
\end{proof}

\vspace{.1in}

\begin{proposition}
For $p,q,r\geq 1$, the elliptic surface $E(p+q+2r+1)$ has a handle decomposition as in Figure~$\ref{fig8}$, where $k=3r-2,3r-1,3r$. The obvious cusp neighborhood is isotopic to the regular neighborhood of a cusp fiber of $E(p+q+2r+1)$. The symbol $F$ denotes the class of a regular fiber of $E(p+q+2r+1)$ in $H_2(E(p+q+2r+1);\mathbf{Z})$.
\begin{figure}[ht!]
\begin{center}
\includegraphics[width=4.3in]{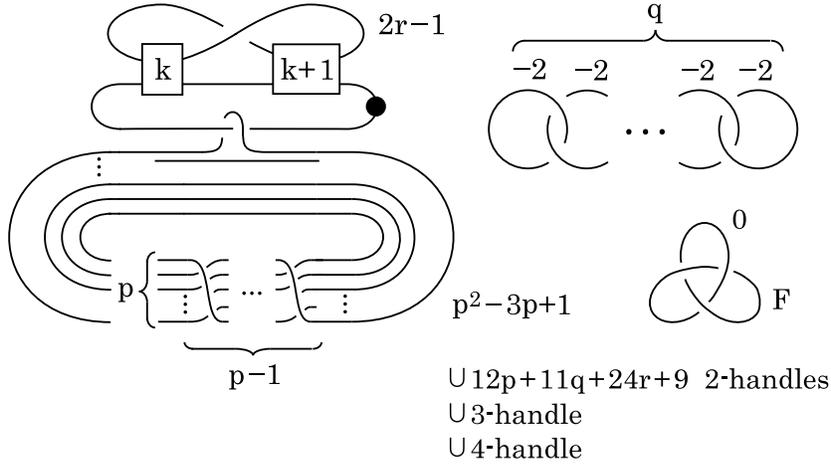}
\caption{$E(p+q+2r+1)$ $(k=3r-2,3r-1,3r)$}
\label{fig8}
\end{center}
\end{figure}
\end{proposition}
\begin{proof}
We give a proof for $k=3r-2$. The other cases are only slightly different. 
Figure~\ref{fig7} shows that $E(p+q+2r+1)$ has a handle decomposition as in Figure~\ref{fig13}. We slide handles as in Figure~\ref{fig14}. Notice that we ignored many 2-handles in the process. We next introduce a 2-handle/3-handle pair and slide the new 2-handle as in the second and third diagrams of Figure~\ref{fig15}. We repeat a handle slide as in the third, fourth and fifth diagrams. (For more details, see \cite[Figure~15--17]{Y1} and \cite[Proposition~3.1.(1)]{Y2}.) We now get the sixth diagram by an isotopy. By sliding the $p^2$-framed knot over $-1$-framed knots, we obtain the last diagram. This procedure gives a diagram of $E(p+q+2r+1)$ as in Figure~$\ref{fig8}$. 
\end{proof}
By blowing ups, we get the following corollary. 

\vspace{.1in}

\begin{corollary}
For $p,q,r\geq 1$, the $4$-manifold $E(p+q+2r+1)\#(2r-1) \overline{\mathbf{C}\mathbf{P}^2}$ has a handle decomposition as in Figure~$\ref{fig9}$. The symbol $F$ denotes the class of a regular fiber of $E(p+q+2r+1)$ in $H_2(E(p+q+2r+1);\mathbf{Z})$. The symbols $E_1,E_2,\dots,E_{2r-1}$ denote a standard basis of $H_2((2r-1) \overline{\mathbf{C}\mathbf{P}^2};\mathbf{Z})$. 
\begin{figure}[ht!]
\begin{center}
\includegraphics[width=4.3in]{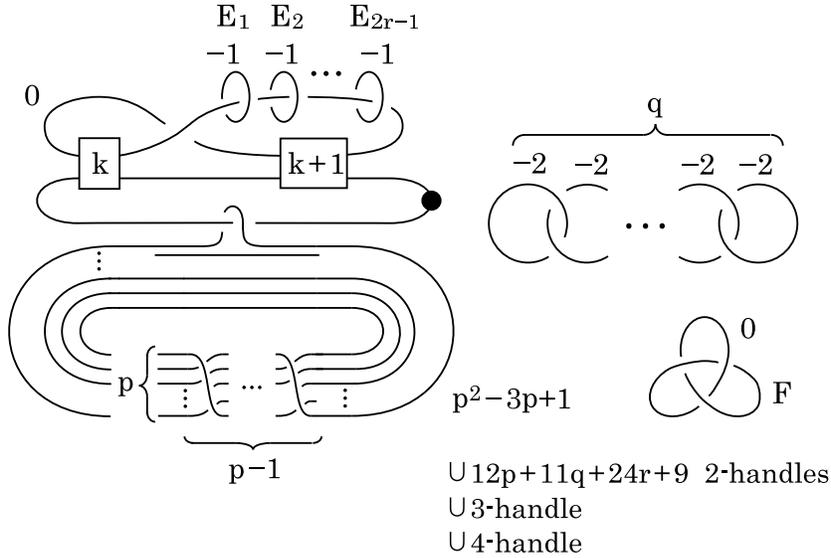}
\caption{$E(p+q+2r+1)\#(2r-1) \overline{\mathbf{C}\mathbf{P}^2}$ \newline $\qquad \qquad \; \; $ $(k=3r-2,3r-1,3r)$}
\label{fig9}
\end{center}
\end{figure}
\end{corollary}
We can now prove Theorem~\ref{th:exotic by cork}. 
\begin{proof}[Proof of Theorem~\ref{th:exotic by cork}]
Let $\alpha$ be a generator of $H_2(C_1(m,n,p,0);\mathbf{Z})\cong \mathbf{Z}$. Let $\beta$ be a generator of $H_2(C_2(m,n,p,0);\mathbf{Z})\cong \mathbf{Z}$. We can easily check that $\beta$ is represented by a surface with its genus $\frac{1}{2}(p^2-3p+2)$. 

\vspace{.05in}

\noindent (1) (i) The case $m=p^2-3p+1$ and $n=3r-2$ $(r\in \mathbf{Z}_{\geq2})$: 
Since $PD((p+2r-1)F)$ is a Seiberg-Witten basic class of $E(p+2r+1)$ (see \cite{GS}.), the blow up formula shows that $K:=PD((p+2r-1)F+E_1+E_2+\dots+E_{2r-1})$ is a Seiberg-Witten basic class of $E(p+2r+1)\#(2r-1) \overline{\mathbf{C}\mathbf{P}^2}$. Here the symbol $PD$ denotes the Poincar\'{e} dual. Let $\alpha'$ be the element of $H_2(E(p+2r+1)\# (2r-1)\overline{\mathbf{C}\mathbf{P}^2};\mathbf{Z})$ induced from $\alpha$ by the inclusion of $C_1(m,n,p,0)$ as shown in Figure~\ref{fig9}. Put $g$ as the genus of a surface which represents $\alpha$. Then the adjunction inequality for the symplectic $4$-manifold $E(p+2r+1)\# (2r-1)\overline{\mathbf{C}\mathbf{P}^2}$ shows the following inequality:
\begin{equation*}
\lvert K(\alpha')\rvert +{\alpha'}^2\leq 
\begin{cases}
2g-2&\text{if\: $g\geq 2$}\\
0&\text{if\: $g=0,1$}.
\end{cases}
\end{equation*}
This inequality together with Figure~\ref{fig9} shows 
\begin{equation*}
g\geq \frac{1}{2}(p^2-3p+2r+2).
\end{equation*}
Note that the right side of the inequality above is greater than or equal to $2$. 
Since $\beta$ is represented by a surface with its genus $\frac{1}{2}(p^2-3p+2)$, the inequality above shows that $C_1(m,n,p,0)$ and $C_2(m,n,p,0)$ are not diffeomorphic. 

(ii) The other cases: We first blow up Figure~\ref{fig9} and  proceed similarly as above by using the adjunction inequality. 
\vspace{.05in}

\noindent (2) The same argument as the proof of (1) shows the claim.
\vspace{.05in}

\noindent (3) Suppose that $C_1(m,n,p,q)$ and $C_2(m,n,p,q)$ are diffeomorphic. Then there exists a genus $\frac{1}{2}(p^2-3p+2)$ surface with self-intersection number $m$ in $C_1(m,n,p,q)$, since $C_2(m,n,p,q)$ has such a surface. Let $\gamma$ be an element of $H_2(C_1(m,n,p,q);\mathbf{Z})$ with self-intersection number $m$. Put $g$ as the genus of a suface which represent $\gamma$. 

(i) The case $m=p^2-3p+1$ and $n=3r-2$ $(r\in \mathbf{Z}_{\geq2})$: Let $\gamma'$ be the element of $H_2(E(p+q+2r+1)\# (2r-1)\overline{\mathbf{C}\mathbf{P}^2};\mathbf{Z})$ induced from $\gamma$ by the inclusion of $C_1(m,n,p,q)$ as shown in Figure~\ref{fig9}. The fact ${\gamma'}^2=m\geq 0$ and Figure~\ref{fig2} shows the existence of some element $x$ of $H_2(E(p+q+2r+1)\# (2r-1)\overline{\mathbf{C}\mathbf{P}^2};\mathbf{Z})$ such that $x\cdot ((p+q+2r-1)F+E_1+E_2+\dots+E_{2r-1})=0$ and $\gamma'=a\alpha'+x$ $(a\geq 1)$. Since $K:=PD((p+q+2r-1)F+E_1+E_2+\dots+E_{2r-1})$ is a Seiberg-Witten basic class of $E(p+q+2r+1)\# (2r-1)\overline{\mathbf{C}\mathbf{P}^2}$, the adjunction inequality gives the following inequality:
\begin{equation*}
\lvert K(\gamma')\rvert +{\gamma'}^2\leq 
\begin{cases}
2g-2&\text{if\: $g\geq 2$}\\
0&\text{if\: $g=0,1$}.
\end{cases}
\end{equation*}
This inequality together with Figure~\ref{fig9} shows 
\begin{equation*}
g\geq \frac{1}{2}(p^2-3p+3+a(2r-1)).
\end{equation*}
Note that the right side of the inequality above is greater than or equal to $2$. 
On the other hand, the assumption shows that $\gamma$ is represented by a genus $\frac{1}{2}(p^2-3p+2)$ surface. This is a contradiction. 

(ii) The other cases: We first blow up Figure~\ref{fig9} and  proceed similarly as above by using the adjunction inequality. 
\end{proof}
\begin{remark}
In Theorem~\ref{th:exotic by cork}, we split the cases into (1) and (2) to avoid the $n\leq 3,\, p=1,2,\, q=0$ case. This is because, in this case, the adjunction inequialities in the proof above cause difficuluties with the genus $g=0$ case. 
\end{remark}
The proof above of Theorem~\ref{th:exotic by cork}.(1) shows the corollary below. This corollary together with Lemma~\ref{lem:basic} clearly gives Theorem~\ref{th:1.4}.
\begin{corollary}If $r\geq 2$, $p\geq 1$ and $m\leq p^2-3p+1$, then 
\begin{equation*}
G(C_1(m,3r-2,p,0))-G(C_2(m,3r-2,p,0))\geq r
\end{equation*}
\end{corollary}
\subsection{Enlarging plugs}
In this subsection, by enlarging plugs $W_{m,n}$, we construct exotic smooth pairs of compact Stein 4-manifolds. We also construct non-homeomorphic pairs of compact Stein 4-manifolds with the same boundary $3$-manifold. 
\begin{definition}
Let $P_1(m,n)$ and $P_2(m,n)$ be the simply connected compact smooth $4$-manifolds in Figure~\ref{fig10}.
\begin{figure}[ht!]
\begin{center}
\includegraphics[width=4.3in]{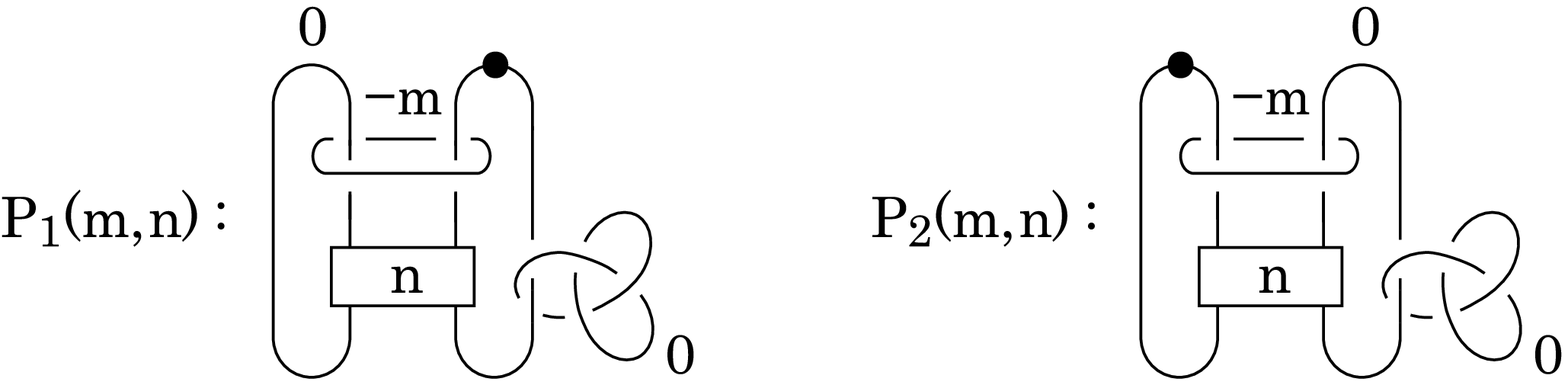}
\caption{}
\label{fig10}
\end{center}
\end{figure}
\end{definition} 
We can easily prove the lemma below, similarly to the proof of Lemma~\ref{lem:basic}. See also the equivalent diagrams of $P_1(m,n)$ and $P_2(m,n)$ in Figure~\ref{fig11}. 
\begin{lemma}
$P_1(m,n)$ and $P_2(m,n)$ are simply connected compact Stein $4$-manifolds for $m,n\geq 1$. The boundaries $\partial P_1(m,n)$ and $\partial P_2(m,n)$ are diffeomorphic. 
\end{lemma}
\begin{figure}[ht!]
\begin{center}
\includegraphics[width=4.3in]{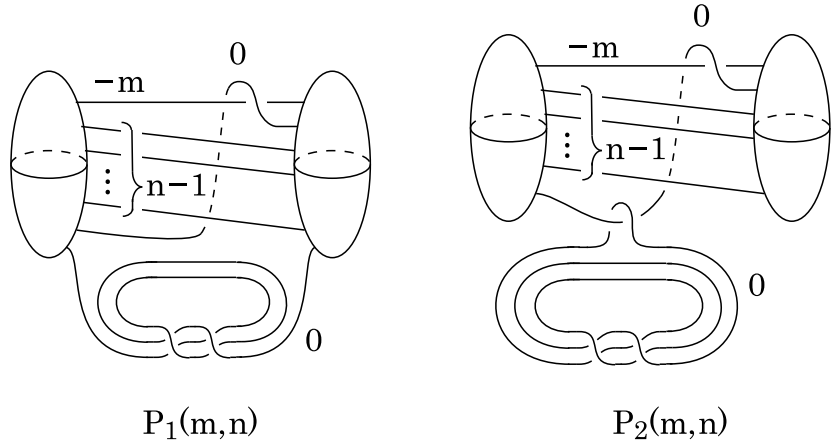}
\caption{$P_1(m,n)$ and $P_2(m,n)$}
\label{fig11}
\end{center}
\end{figure}
We first prove the following theorem. 
\begin{theorem}\label{th:exotic plug}
Simply connected compact Stein $4$-manifolds $P_1(1,3)$ and $P_2(1,3)$ are homeomorphic but not diffeomorphic to each other.
\end{theorem}
\begin{proposition}
$E(2)\# 2\overline{\mathbf{C}\mathbf{P}^2}$ has a handle decomposition as in Figure~\ref{fig12}. Here $E_1$ and $E_2$ denote a standard basis of $H_2(2\overline{\mathbf{C}\mathbf{P}^2};\mathbf{Z})$. 
\begin{figure}[ht!]
\begin{center}
\includegraphics[width=2.5in]{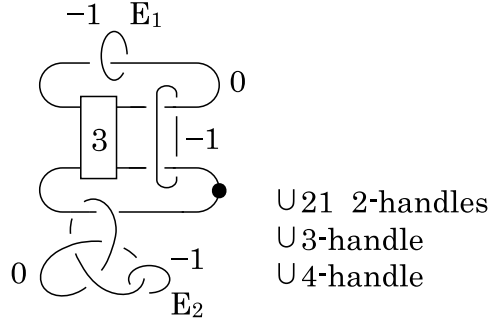}
\caption{$E(2)\# 2\overline{\mathbf{C}\mathbf{P}^2}$}
\label{fig12}
\end{center}
\end{figure}
\end{proposition}
\begin{proof}
The $n=2$ case of Figure~\ref{fig7} obviously gives the first diagram of $E(2)$ in Figure~\ref{fig17}. We get the second diagram of Figure~\ref{fig17} by handle moves as in the $p=3$ case of Figure~\ref{fig15}. Handle slides as in Figure~\ref{fig16} gives the third diagram of Figure~\ref{fig17}. Finally blow up twice. Then we get Figure~\ref{fig12}.
\end{proof}
\begin{proof}[Proof of Theorem~$\ref{th:exotic plug}$]
We can easily check that the intersection forms of both $P_1(1,3)$ and $P_2(1,3)$ are $\langle 1 \rangle\oplus \langle -1 \rangle$. Boyer's theorem~\cite{B} thus shows that $P_1(1,3)$ and $P_2(1,3)$ are homeomorphic. 

It follows from Figure~\ref{fig12} that the subspace $H_2(P_1(1,3);\mathbf{Z})$ of $H_2(E(2)\#2\overline{\mathbf{C}\mathbf{P}^2};\mathbf{Z})$ is spanned by elements $x_1$ and $x_2$ such that $E_1\cdot x_1=1$, $E_2\cdot x_1=0$, $E_1\cdot x_2=0$ and $E_2\cdot x_2=1$. The blow-up formula shows that the Seiberg-Witten basic classes of $E(2)\# 2\overline{\mathbf{C}\mathbf{P}^2}$ are $\pm PD(E_1\pm E_2)$. Using the adjunction inequality for $E(2)\# 2\overline{\mathbf{C}\mathbf{P}^2}$, we can now easily prove that there exists no non-zero element of $H_2(P_1(1,3);\mathbf{Z})$ such that the square of the element is zero and that the element is represented by a torus.

On the other hand, we can easily check that there exists a non-zero element of $H_2(P_2(1,3);\mathbf{Z})$ such that the square of the element is zero and that the element is represented by a torus. Therefore, the claim follows. 
\end{proof}
\begin{remark}
We can easily extend Theorem~\ref{th:exotic plug} as in Theorem~\ref{th:exotic by cork}. 
\end{remark}
We finally give the theorem below. This shows Theorem~\ref{th:1.3}. 
\begin{theorem}
If $m\geq 1$ is odd, and $n\geq 1$ is even, then simply connected compact Stein $4$-manifolds $P_1(m,n)$ and $P_2(m,n)$ have non-isomorphic intersection forms, whereas they have diffeomorphic boundaries and isomorphic integer homology groups.
\end{theorem}
\begin{proof}
We can easily check that the intersection form of $P_1(m,n)$ is odd and that the one of $P_2(m,n)$ is even, under the assumption. 
\end{proof}

\begin{figure}[p]
\begin{center}
\includegraphics[width=4.5in]{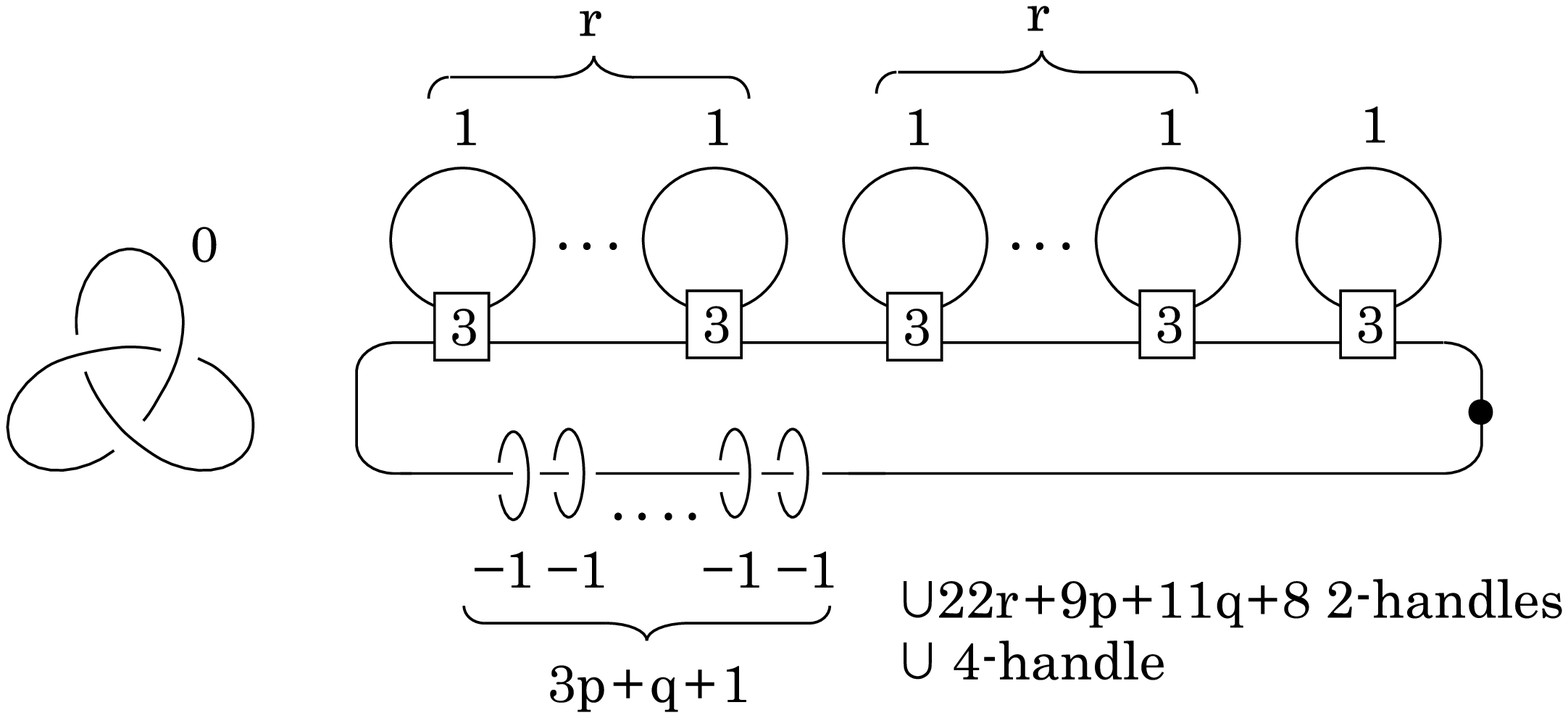}
\caption{$E(p+q+2r+1)$}
\label{fig13}
\end{center}
\end{figure}
\begin{figure}[]
\begin{center}
\includegraphics[width=4.9in]{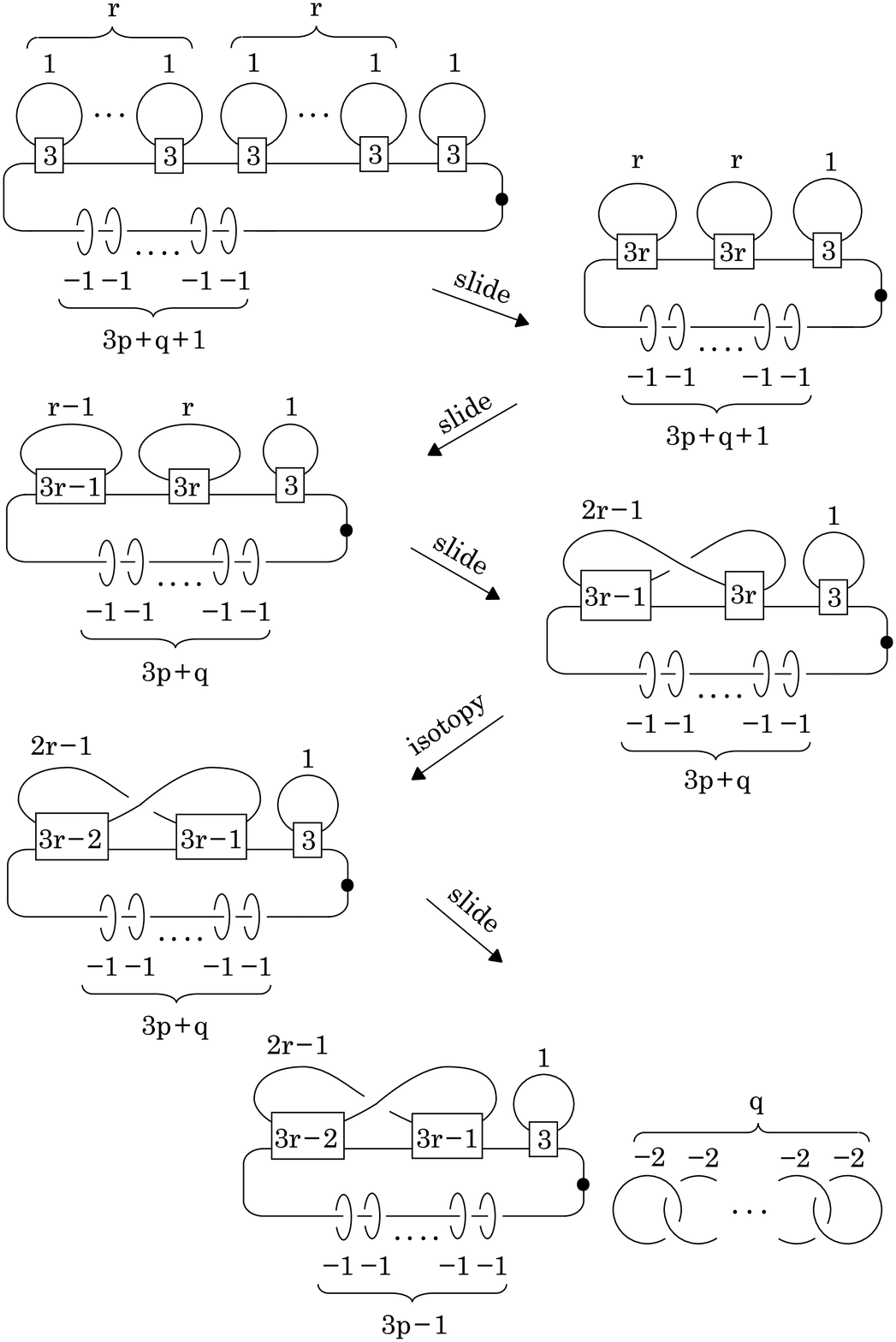}
\caption{handle slides}
\label{fig14}
\end{center}
\end{figure}
\begin{figure}[]
\begin{center}
\includegraphics[width=4.9in]{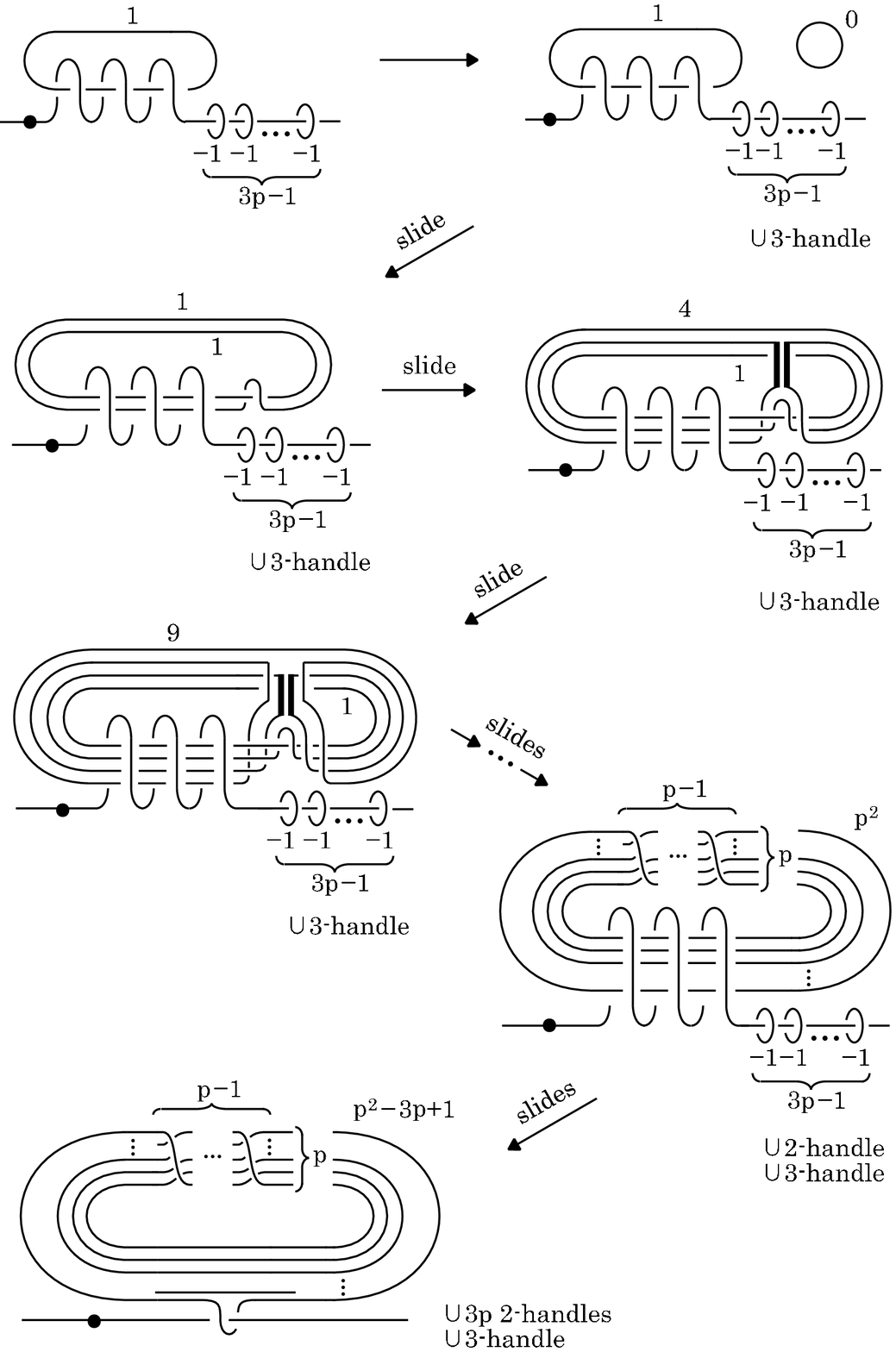}
\caption{handle slide}
\label{fig15}
\end{center}
\end{figure}
\begin{figure}[]
\begin{center}
\includegraphics[width=4.9in]{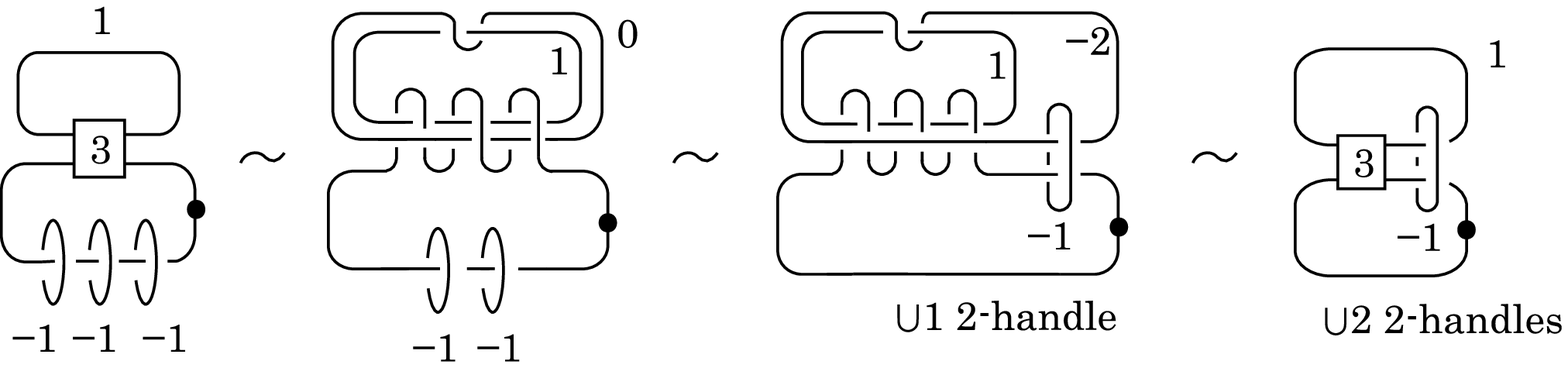}
\caption{handle slide}
\label{fig16}
\end{center}
\end{figure}
\begin{figure}[]
\begin{center}
\includegraphics[width=4in]{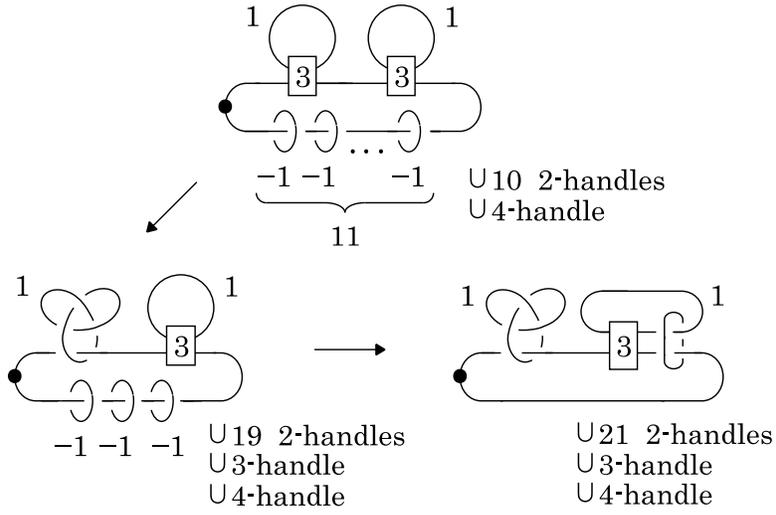}
\caption{handle moves of $E(2)$}
\label{fig17}
\end{center}
\end{figure}

\end{document}